\documentclass{article}
\usepackage[utf8]{inputenc}
\usepackage[english]{babel}
\usepackage{booktabs}
\usepackage{mathrsfs}
\usepackage{mathtools}
\usepackage{amsmath}
\usepackage{amssymb}
\usepackage{graphicx}
\usepackage{amsfonts}
\usepackage{amsthm}
\usepackage{mathrsfs}

\usepackage{color}

\usepackage[mathscr]{euscript}
\usepackage{subfiles}
\usepackage{enumitem}

\usepackage{centernot}
\usepackage{bm}
\usepackage{array}
\usepackage{bbm}
\usepackage{authblk}

\newcommand{\R}{\mathbb{R}} 
\newcommand{\D}{\mathbb{D}} 

\newcommand{\N}{\mathbb{N}} 
\newcommand{\Fe}{\mathbb{F}}

\renewcommand{\tilde}[1]{\widetilde{#1}} 

\renewcommand{\P}{\mathbb{P}} 

\newtheorem{theorem}{Theorem}\newtheorem{corollary}{Corollary}\newtheorem{lemma}{Lemma}\theoremstyle{definition}
\theoremstyle{remark}

\usepackage{tikz}
\usepackage{pgfplots}

\newcommand{\PFspace}[1]{\left(\Omega, \mathcal{F}, \left(\mathcal{F}_t\right)_{#1}, \mathbb{P}\right)} 

\usepackage{subcaption}

\numberwithin{equation}{section}

\begin{document}
	\title{Estimating the Copula of a class of Time-Changed Brownian Motions:
		A non-parametric Approach}
	
	\author{Orimar Sauri\thanks{osauri@math.aau.dk} }
	\author{Toke C. Zinn\thanks{tokecz@math.aau.dk}}
	\affil{Department of Mathematical Sciences, Aalborg University \protect\\
		Skjernvej 4A, 9220, Aalborg, Denmark} 
	\date{This Version: \today}
	
	\maketitle
\begin{abstract}
	Within a high-frequency framework, we propose a non-parametric
	approach to estimate a family of copulas associated to a time-changed
	Brownian motion. We show that our estimator is consistent
	and asymptotically mixed-Gaussian. Furthermore, we test its finite-sample
	accuracy via Monte Carlo.
	
\end{abstract}

\section{Introduction}
One of the most fundamental results in probability theory is the so-called Sklar's Theorem. It states that for every random vector $\mathbf{X}=(X_{1},\ldots,X_{n})$
there exists a copula $C$ (see Section 2) such that 
\begin{align*}
H(x)=C(F(x)),\quad\forall x\in\R^{n},
\end{align*}
where $H:\R^{n}\to[0,1]$ is the cumulative distribution function
(cdf for short) of $\mathbf{X}$ and
\begin{align}
F(x):=(F_{1}(x_{1}),F_{2}(x_{2}),\dots,F_{n}(x_{n})),\quad x\in\R^{n},
\end{align}
and $F_i$ is the cdf of $X_i$, for $i$ in $\{1,2,\dots,n\}$.

Since then, copulas have been applied in a large number of sciences, for instance in hydrology (\cite{salvadori2007use}), engineering (\cite{yang2017reliability}), and, perhaps most noticeably, in finance.  In finance, their primary use is in risk management and portfolio allocation. Specifically, copulas are used to model the joint distribution of financial assets. For a detailed account see \cite{Cherubinietal12} and \cite{genest2009advent}. Note that this way of using copulas can be considered as a spatial way of modeling dependence, i.e. describing the dependence between two or several stochastic processes at distinct times. Copulas have also found their use in temporal modelling as well. \cite{darsow1992copulas} characterized Markov processes by means of their copulas. \cite{wilson2010copula} also define so-called copula processes. Furthermore, \cite{darsow1992copulas} derived the copula of the bivariate distributions of a Brownian motion. In addition, \cite{schmitz2003copulas}, argued that the copula of a Brownian motion could be uses to construct similar processes with arbitrary marginal distributions. 

In the spatial set-up, parametrical and non-parametrical inference for copulas is well documented. See for instance \cite{scaillet2002nonparametric} and \cite{genest1995semiparametric}. However, in the temporal case very little statistical analysis has been done, see for instance \cite{chen2006estimation}. This paper aims at developing some results in that direction.

In the present work, we concentrate on the statistical inference for a family of copulas associated with the finite-dimensional distributions of a class of time-changed Brownian motions. More precisely, we propose a non-parametric estimator for a family of conditional copulas linked to a time-changed Brownian motion. We show consistency and asymptotic (mixed) normality under the assumption that the process is observed in a high-frequency set-up.

The paper is structured as follows. Section \ref{sec:Background} introduces the notations used through the paper and discusses some essential preliminaries. In Section \ref{sec:Results} we present our results and we show the performance of the estimator in finite-samples. In the last section the proofs of our main results are presented.

\section{Background\label{sec:Background}}

In this section we recall several definitions and properties required
to present our main results. Throughout this paper $\PFspace{[0,\infty)}$
will denote a filtered probability space satisfying the usual assumptions
of right-continuity and completeness. If $\mathcal{X}$ is a set,
then we will denote the $n$-fold cartesian product of $\mathcal{X}$
with itself as $\mathcal{X}^{n}$. Similarly, if $(\mathcal{X},d)$
is a metric space, then the product metric defined on $\mathcal{X}^{n}$
is denoted $d^{n}$.

\subsection{Copulas of Time-Changed Brownian motion}

For $a =(a_1,\dots,a_n),b = (b_1,\dots,n)\in\R^{n}, n\in\N$ we write $a\leq b$ ($a<b$) if $a_{i}\leq b_{i}$ ($a_{i}<b_{i}$) for every $i\in\{1,\dots,n\}$. Whenever $a\leq b$ with $a,b\in\mathbb{R}^n$, we
define $(a,b]:=\bigtimes_{i=1}^{n}(a_{i},b_{i}]$. Recall that a\textit{
copula} is a function $C:[0,1]^{n}\to[0,1]$ satisfying the
following:
\begin{equation}
\begin{aligned}C(u) & =0,\quad u\in[0,1]^{n}:\exists i\in\{1,2,\dots,n\}:u_{i}=0\\
	C(1,1,\dots,u,\dots,1) & =u,\quad u\in[0,1]\\
	\int_{(a,b]}dC & \geq0,\quad\forall a,b\in[0,1]^{n}:a\leq b.
\end{aligned}
\label{copuladef}
\end{equation}
In this paper we focus on the time-changed Brownian motion 
\begin{align}
X_{t}=W_{T_{t}},\quad t\geq0,\label{eq:CLM}
\end{align}
where $(T_{t})_{t\geq0}$ is a continuous random time change independent
of $W$, that is, it is a non-decreasing process taking values in
$[0,\infty]$ such that $T_{t}$ is a $(\mathcal{F}_{t})_{t\geq0}$-stopping
time for all $t\geq0$. Within this framework, for all $0\leq t_{0}<t_{1}<t_{2}<\cdots<t_{n}$
the copula associated to $(X_{t_{0}},X_{t_{1}},\ldots,X_{t_{n}})$,
is completely determined by the law of $\mathcal{T}_{n}:=(T_{t_{0}},\ldots,T_{t_{n}})$
and the random fields 
\begin{align}
C(t_{k-1},t_{k};u,v):=\psi(T_{k-1},T_{k};u,v),\,\,u,v\in[0,1],k=1,\ldots,n,\label{eq:CopulaTCBM}
\end{align}
where 
\begin{equation}
\psi(s,t;u,v):=\begin{cases}
	\int_{0}^{u}\Phi\left(\frac{\sqrt{t\lor s}\Phi^{-1}(v)-\sqrt{s\land t}\Phi^{-1}(w)}{\sqrt{\left|t-s\right|}}\right)dw, & \text{if }\left|t-s\right|>0;\\
	u\land v & \text{if }t=s>0;\\
	uv & \text{otherwise},
\end{cases}\label{bmcop}
\end{equation}
in which $\Phi$ denotes the cdf of the standard normal distribution.
For more details on the previous statements we refer the reader to
\cite{Cherubinietal12} and \cite{darsow1992copulas}. It is not difficult
to see that the mapping $(u,v)\mapsto C(t,s;u,v)$ satisfies (\ref{copuladef})
almost surely for all $t\geq s$. Moreover, it holds that 
\begin{align}\label{eq:condcopspecification}
C(s,t;F_{s}(x\mid T_{t},T_{s}),F_{t}(y\mid T_{t},T_{s}))=\mathbb{P}(X_{s}\leq x,X_{t}\leq y\mid T_{t},T_{s}),
\end{align}
where $F_{t}(x\mid T_{t},T_{s})$ denotes the cdf of $X_{t}$ given
$(T_{t},T_{s})$. Motivated by \eqref{eq:condcopspecification} and the terminology used in \cite{patton2006modelling},
we will refer to $C(t,s;\cdot,\cdot)$ as \textit{the conditional
copula associated to $X$}.
\subsection{Limit Theorems and Convergence}

The notations $\overset{\mathbb{P}}{\rightarrow}$ and $\overset{d}{\rightarrow}$
stand, respectively, for convergence in probability and in distribution
of random vectors. As usual the space of cádlág fields will be denoted by $\D([0,\mathscr{T}]^{n};\mathbb{R}^{d})$. If $X$ and $X^{n}$ are two cádlág processes we
write $X^{n}\overset{u.c.p.}{\Rightarrow}X$, whenever 
\begin{align*}
\lim_{n\to\infty}\P(\sup_{0\leq t\leq \mathcal{\mathscr{T}}}\left\Vert X_{t}^{n}-X_{t}\right\Vert \geq\varepsilon)=0,\quad\forall\,\mathcal{\mathscr{T}},\varepsilon>0.
\end{align*}
A sequence of random vectors $(\xi_{n})_{n\geq1}$ on $\left(\Omega,\mathcal{F},\mathbb{P}\right)$
is said to\textit{ converges stably in law} towards $\xi$ (in symbols
$\xi_{n}\overset{s.d}{\longrightarrow}\xi$), which is defined on
an extension of $\left(\Omega,\mathcal{F},\mathbb{P}\right)$, say
$\left(\tilde{\Omega},\tilde{\mathcal{F}},\tilde{\mathbb{P}}\right)$,
if for every continuous and bounded function $f$ and any bounded
random variable $\chi$ it holds that 
\[
\mathbb{E}(f(\xi_{n})\chi)\rightarrow\tilde{\mathbb{E}}(f(\xi)\chi),
\]
where $\tilde{\mathbb{E}}$ denoted expectation w.r.t. $\tilde{\mathbb{P}}$.
For a concise exposition of stable convergence see \cite{hausler2015stable}.
Given a stochastic process $Z=(Z_{t})_{t\geq0}$, we will use the
notation $\Delta_{i}^{n}Z:=Z_{i/n}-Z_{(i-1)/n}$, $i\in\N$. The \textit{realized variation} of a process $Z=(Z_{t})_{t\geq0}$
is defined and denoted as the process 
\begin{align*}
[Z]_{t}^n=\sum_{i=1}^{[nt]}(\Delta_{i}^{n}Z)^{2},\quad t\geq0.
\end{align*}
where $[x]$ denotes the integer part of $x\in\R$. It is well known that
if $Z$ is a continuous semimartingale, then $[Z]^n\overset{u.c.p}{\Rightarrow}[Z].$

\section{Estimating the conditional copula of $X$\label{sec:Results}}

As discussed in Section \ref{sec:Background}, the copula associated
to $(X_{t_{0}},X_{t_{2}},\ldots,X_{t_{n}})$, for $0<t_{1}<t_{2}<\cdots<t_{n}$,
is completely determined by the law of $(T_{t_{0}},\ldots,T_{t_{n}})$
as well as the family of conditional copulas $\{C(s,t;u,v):0\leq s,t\leq T,u,v\in[0,1]\}$,
where $C$ as in (\ref{eq:CopulaTCBM}). For the rest of this section,
we propose a non-parametric approach for estimating the later. Our
sample scheme is as follows: The process $X$ is observed on a fixed
interval $[0,\mathscr{T}]$, $\mathscr{T}>0$, at times $t_{i}=i/n$,
for $i=0,1,\ldots,[n\mathscr{T}]$. Thus, motivated by (\ref{eq:CopulaTCBM})
and the fact that $[X]^n\overset{u.c.p}{\Rightarrow}[X]=T$, as
$n\rightarrow\infty$, we propose to estimate $C$ via 
\begin{align*}
C^{n}(s,t;u,v):=\psi([X]_{s}^n,[X]_{t}^n;u,v),\quad u,v\in[0,1],0\leq s,t.
\end{align*}
Our first result shows that $C^{n}$ is indeed a consistent estimator
for $C$.

\begin{theorem}\label{consistencythm}Let $t_{0}\geq0$, such that
$\mathbb{P}(T_{t_{0}}>0)=1$. Then for all $u,v\in[0,1]$, $\mathscr{T}>t_{0}$, and $\varepsilon>0$, it holds

\begin{align*}
	\P\left(\sup_{t_{0}\leq s,t\leq\mathscr{T}}|C(s,t;u,v)-C^{n}(s,t;u,v)|\geq\varepsilon\right)\rightarrow0,\,\,\,n\rightarrow\infty.
\end{align*}
\end{theorem}

Now we proceed to derive second-order asymptotics for $C^{n}$. In
order to do this, we require stronger assumptions on the structure
of $X$. Specifically, we are going to assume that there is a $(\mathcal{G}_{t})_{t\geq0}$-Brownian
motion $B$ such that 
\begin{equation}
X_{t}=\int_{0}^{t}\sigma_{s}dB_{s},\,\,\,t\geq0,\label{dynamicXCLT}
\end{equation}
where $\sigma$ is a cádlág process. Observe that by Knight's Theorem,
$X$ admits the representation 
\begin{equation}
X_{t}=W_{T_{t}},\,\,\,t\geq0,\label{timechanged}
\end{equation}
where $T_{t}=\int_{0}^{t}\sigma_{r}^{2}dr,\,\,t\geq0$. Thus, if $\sigma$ is assumed to be independent of $B$, we can further
choose $W$ to be independent of $\sigma$. This can be seen easily
in the case when $\int_{0}^{\infty}\sigma_{r}^{2}dr=\infty$. Indeed,
in that situation it is well known (see for instance \cite{BNShirChangetimebook})
that (\ref{timechanged}) holds with 
\[
W_{t}=X_{A_{t}},\,\,\,t\geq0,
\]
where $A_{t}=\inf\{s\geq0:T_{s}>t\}$, which easily implies that $W$ is independent of $\sigma$. The general
case can be analysed in a similar way. Assuming that $X$ admits the
representation (\ref{dynamicXCLT}) we obtain the following Central
Limit Theorem for $C^{n}.$

\begin{theorem}\label{clthm}Let $X$ be given by (\ref{dynamicXCLT})
and assume that $\mathbb{P}(\sigma_{t}^{2}>0)=1$ for all $t\geq0$.
Fix $u,v\in(0,1)$ and denote by $\nabla\psi(s,t;u,v)=(\partial_{t}\psi(s,t;u,v),\partial_{s}\psi(s,t;u,v))$.
Then for $0<s<t$, as $n\rightarrow\infty$ 
\[
\sqrt{n}\left[C^{n}(s,t;u,v)-C(s,t;u,v)\right]\xrightarrow{s.d}\sqrt{V_{s,t}}N(0,1),
\]
where $N(0,1)$ is a standard normal random variable independent of
$\mathcal{F}$ and 
\[
V_{s,t}=2\nabla\psi(T_{s},T_{t};u,v)\begin{bmatrix}Q_{t} & Q_{s}\\
	Q_{s} & Q_{s}
\end{bmatrix}\nabla\psi(T_{s},T_{t};u,v)^{\prime},
\]
in which $Q_{t}:=\int_{0}^{t}\sigma_{r}^{4}dr$, $t\geq0.$
\end{theorem}

A simple way to estimate $V_{s,t}$ is by using power variations:
If $\sigma$ is càdlàg (see for instance \cite{jacod2011discretization}),
then as $n\rightarrow\infty$
\[
Q^n_{t}:=\frac{n}{3}\sum_{i=1}^{[nt]}\left|\Delta_{i}^{n}Z\right|^{4}\overset{u.c.p}{\Rightarrow}Q_{t}.
\]
Thus, a feasible estimator for $V_{s,t}$ is 
\[
V^n_{s,t}:=2\nabla\psi([X]^n_{s},[X]^n_t;u,v)\begin{bmatrix}Q^n_{t} & Q^n_{s}\\
Q^n_{s} & Q^n_{s}
\end{bmatrix}\nabla\psi([X]^n_{s},[X]^n_t;u,v)^{\prime}\overset{\P}{\rightarrow} V_{s,t}.
\]
Thus, we have an easy consequence of the previous theorem:

\begin{corollary}\label{corollaryclt}Under the assumptions of the
previous theorem we have that 
\[
\sqrt{\frac{n}{V^n_{s,t}}}\left[C^{n}(s,t;u,v)-C(s,t;u,v)\right]\xrightarrow{s.d}N(0,1).
\]

\end{corollary}

\subsection{Simulation study}

In this part we study the finite-sample behavior of our proposed estimator.
We use here Monte Carlo simulations to investigate the sensitivity
of $C^{n}$ to the variation of $s,t,u,v$ as well as the sample size.
Our set-up is as follows: The volatility term $\sigma^{2}$ is simulated
according to the so-called Cox-Ingersoll-Ross process, i.e. $\sigma^{2}$ satisfies
the stochastic differential equation 
\begin{align*}
d\sigma_{t}^{2}=\kappa(\theta-\sigma_{t}^{2})dt+\nu\sqrt{\sigma_{t}^{2}}dW_{t},\quad\sigma_{0}^{2}=s_{0}.
\end{align*}
The parameters are set $(\kappa,\theta,\nu,s_{0})=(0.5,1.5,1,1.5)$
in such a way that the Feller condition $2\kappa\theta>\nu^{2}$ is satisfied.
Based on this, we sample over the interval $[0,1]$ equidistant discretizations
$(X_{i/n})_{i=1}^{n}$, where $X$ is given as in (\ref{timechanged}). 

\begin{figure}[htbp]
\begin{subfigure}[b]{.5\textwidth} \includegraphics[scale=0.45]{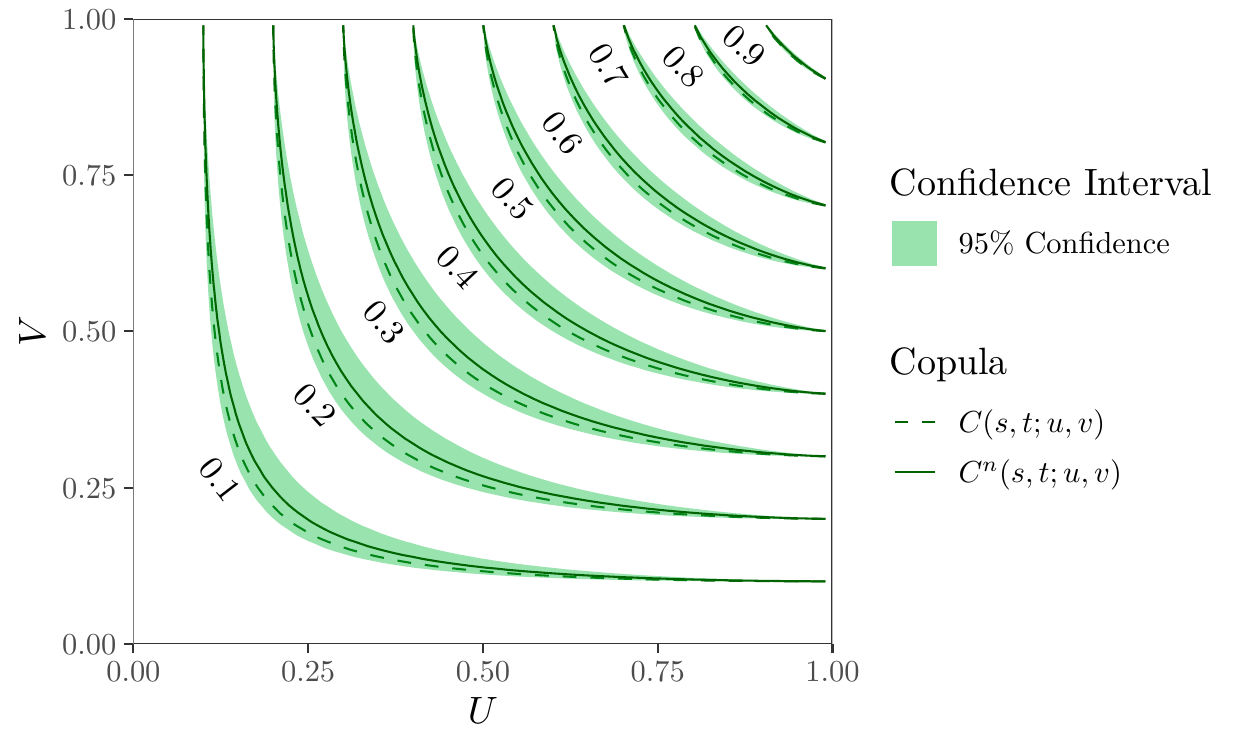}
	\caption{Level curves with estimated confidence intervals for $n=100$ points
		per path.}
\end{subfigure} \hfill{}\begin{subfigure}[b]{.5\textwidth}
	\centering \includegraphics[scale=0.45]{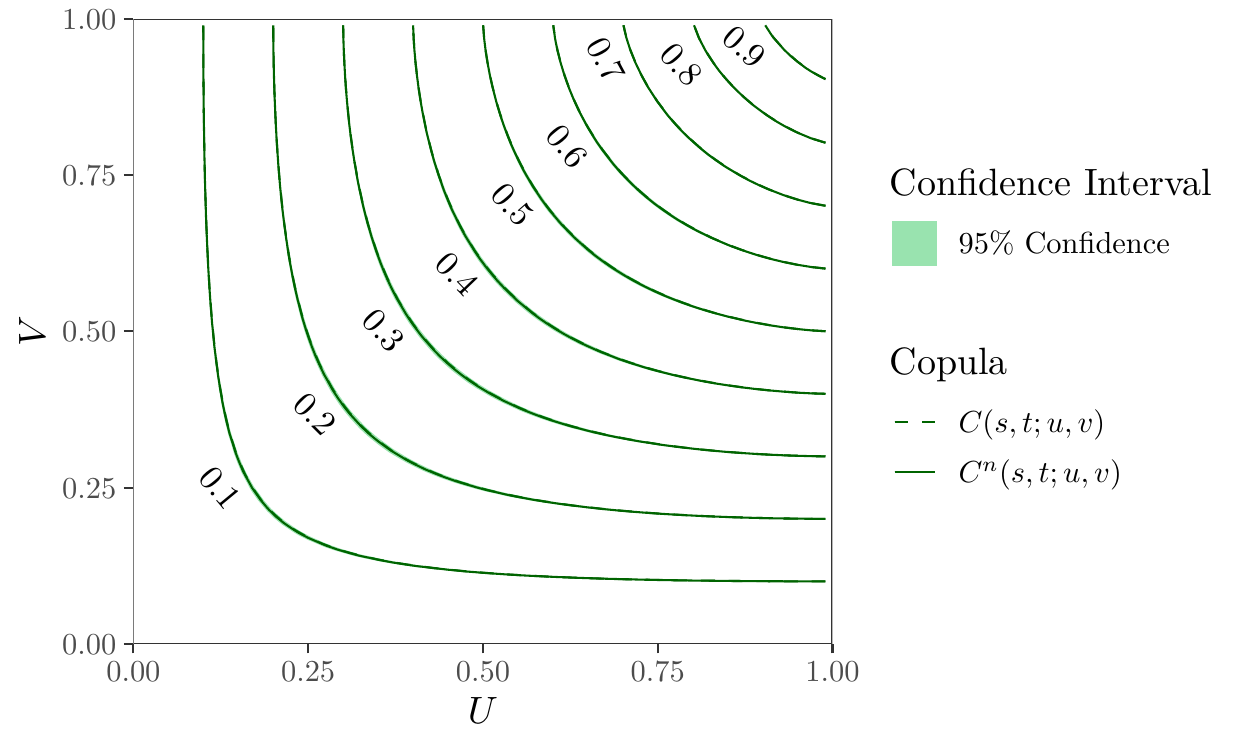} \caption{Level curves with estimated confidence intervals for $n=10000$ points
		per path.}
\end{subfigure} \caption{Level curves for the estimated copula with confidence intervals ($s=0.3$
	and $t=0.7$).}
\label{fig:contours}
\end{figure}

In Figure \ref{fig:contours} we have plotted the level sets of $C^n(s,t;\cdot,\cdot)$ and $C(s,t;\cdot,\cdot)$ together with the $95\%$ confidence contours provided by Theorem \ref{clthm}. The confidence countours behave as we would imagine; near the line $u = v$ we see that that they allow for the largest deviations. As $(u,v)$ approach the boundary of $[0,1]^2$, the intervals diminish as expected: $C^n(s,t;u,v)$ and $C(s,t;u,v)$ coincide on the boundary. Furthermore, as $n$ increases we see that the confidence contours become increasingly narrow as expected. 

\begin{figure}[htbp]
\centering \includegraphics[scale=0.6]{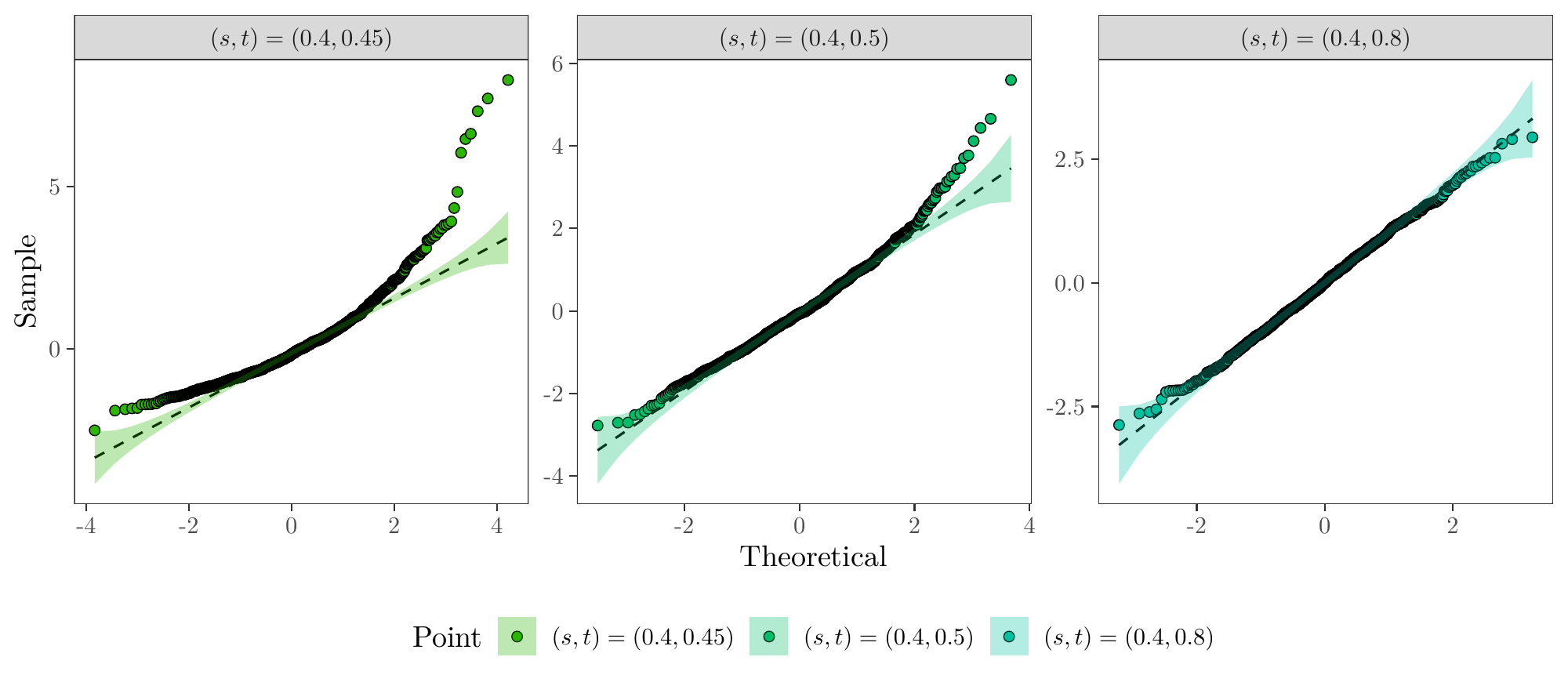} \caption{QQ plot for $\sqrt{n/V^n_{s,t}}\cdot(C^{n}(s,t;u,v)-C(s,t;u,v))$,
	$(n=10000,(u,v)=(0.7,0.3))$.}
\label{fig:QQ}
\end{figure}

In Figure \ref{fig:QQ} we report the finite-sample distribution of our standardized error against the standard normal distribution. We can see that the accuracy of our statistic is sensitive to the boundary points where $s$ and $t$ are close. Again, this behaviour is not unexpected; the temporal gradient $\nabla \psi$, as in Theorem \ref{clthm}, is given by
\[  \begin{bmatrix}\int_0^u\varphi\left(\frac{\sqrt{t}\Phi^{-1}(v)-\sqrt{s}\Phi^{-1}(w)}{\sqrt{t-s}}\right)\left(\frac{\Phi^{-1}(v)}{2\sqrt{t(t-s)}}-\frac{\sqrt{t}\Phi^{-1}(v)-\sqrt{s}\Phi^{-1}(w)}{2\sqrt{t-s}^{3}}\right)dw\\
\\
\int_0^u\varphi\left(\frac{\sqrt{t}\Phi^{-1}(v)-\sqrt{s}\Phi^{-1}(w)}{\sqrt{t-s}}\right)\left(\frac{\sqrt{t}\Phi^{-1}(v)-\sqrt{s}\Phi^{-1}(w)}{2\sqrt{t-s}^{3}}-\frac{\Phi^{-1}(w)}{2\sqrt{s(t-s)}}\right)dw
\end{bmatrix}. \]
Here, $\varphi$ is the density function of a standard Gaussian. We see, that terms proportional to $1/\sqrt{t-s}$ appear. However, recall that as $(s,t) \to (t_0,t_0)$ for $t_0 > 0$ the copula reduces to $C(t_0,t_0;u,v) = u\land v$. Furthermore, it is also very likely that numerical errors influence the result here, due to terms such as $1/\sqrt{t-s}$.

We conclude this section by investigating whether the convergence in Theorem \ref{consistencythm} can be extended to uniform convergence over
$(u,v)$. Specifically, we investigate, via Monte Carlo simulations,
the asymptotic behavior of the statistic
\begin{align*}
\rho(C^{n},C):=\sup_{\substack{(u,v)\in[0,1]^{2}\\
		\tau\leq s,t\leq\mathscr{T}
	}
}|C^{n}(s,t;u,v)-C(s,t;u,v)|,
\end{align*}
as the sample size increases.

\begin{figure}[htbp]
\centering
\input{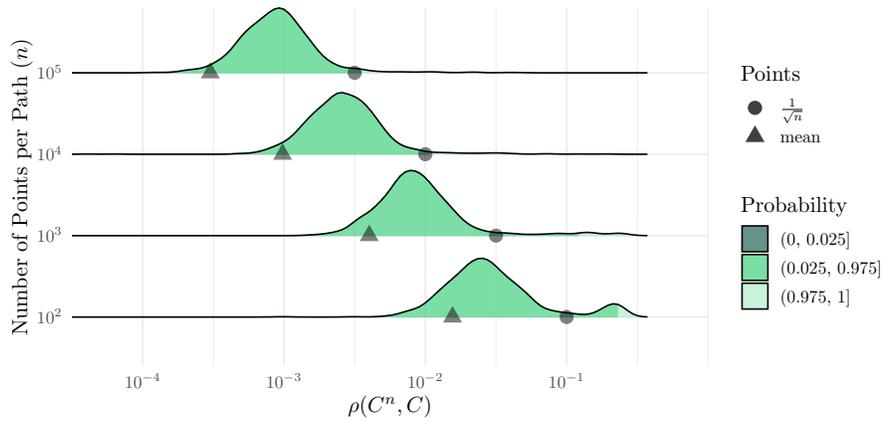} 
\caption{Density estimates for $\rho(C^{n},C)$ on a logarithmic scale.}\label{fig:rhostat}
\end{figure}

Figure \ref{fig:rhostat} shows the density estimates, obtained via a Gaussian kernel density estimate, for $\rho(C^n,C)$ on a logarithmic scale focusing on the peaks. Observe that the mean is located in the lower tail, due to a substantial number of simulations resulting in $\rho(C^n,C)$ being very close $0$. Similarly, $1/\sqrt{n}$ is also added to the plot, showing that the shift in the distribution is, relatively, proportional to $1/\sqrt{n}$ on a logarithmic scale. Based on Figure \ref{fig:rhostat}, these simulations indicate that Theorem \ref{consistencythm} may be extended to include the supremum over $(u,v)\in[0,1]^2$.

\section{Proofs}
The following lemmas are key for the proof of our main results.
\begin{lemma}\label{lemma:limits} Let $\psi$ be as in (\ref{bmcop}).
For all $u,v\in[0,1]$, the mapping $(s,t)\mapsto\psi(s,t;u,v)$ is
continuous in $(0,\infty)^{2}$ and continuously differentiable in
$\{(s,t)\mid 0<s<t\}$.\end{lemma} 

\begin{proof} If $(u,v)$ is in the boundary of $[0,1]^{2}$ the
result is trivial. Suppose that $u,v\in(0,1)$ and put 
\[
g(t,s,w):=\Phi\left(\frac{\sqrt{t\lor s}\Phi^{-1}(v)-\sqrt{t\land s}\Phi^{-1}(w)}{\sqrt{\left|t-s\right|}}\right),\,\,\,t,s>0,w\neq v,
\]
and $g(t,s,w)\equiv0$ when $w=v$. From Example 5.32 in \cite{schmitz2003copulas},
it follows that for almost all $w\in[0,u]$
\begin{align}
	\lim_{(t,s)\rightarrow(t_{0},s_{0})}g(t,s,w)=\begin{cases}
		\mathbbm{1}_{[0,v]}(w) & \text{if }t_{0}=s_{0};\\
		g(t_{0},s_{0},w) & \text{if }t_{0}\neq s_{0}.
	\end{cases}\label{eq:AlmostIndicator}
\end{align}
The continuity then follows by the Lebesgue's Dominated Convergence
Theorem. On the other hand, for $w\neq v$ we have that for $0<s<t$
\[
\begin{bmatrix}\dfrac{\partial g(s,t,w)}{\partial t}\\
	\\
	\dfrac{\partial g(s,t,w)}{\partial s}
\end{bmatrix}=\begin{bmatrix}\varphi\left(\frac{\sqrt{t}\Phi^{-1}(v)-\sqrt{s}\Phi^{-1}(w)}{\sqrt{t-s}}\right)\left(\frac{\Phi^{-1}(v)}{2\sqrt{t(t-s)}}-\frac{\sqrt{t}\Phi^{-1}(v)-\sqrt{s}\Phi^{-1}(w)}{2\sqrt{t-s}^{3}}\right)\\
	\\
	\varphi\left(\frac{\sqrt{t}\Phi^{-1}(v)-\sqrt{s}\Phi^{-1}(w)}{\sqrt{t-s}}\right)\left(\frac{\sqrt{t}\Phi^{-1}(v)-\sqrt{s}\Phi^{-1}(w)}{2\sqrt{t-s}^{3}}-\frac{\Phi^{-1}(w)}{2\sqrt{s(t-s)}}\right)
\end{bmatrix},
\]
where $\varphi$ is density of a standard normal distribution. Since
for any constants $c,k$ such that $k\neq0$, it holds that $\varphi(c+kx)x\rightarrow0$
as $\left|x\right|\rightarrow\infty$, we deduce that
\begin{equation}
	\sup_{w\in[0,1]}\left\Vert \left(\dfrac{\partial g(s,t,w)}{\partial t},\dfrac{\partial g(s,t,w)}{\partial s}\right)\right\Vert <\infty,\,\,\,0<s<t.\label{suppartial}
\end{equation}
Interchanging roles between $s$ and $t$ allow us to conclude that
(\ref{suppartial}) is fulfilled for all $(s,t)\in \{(s,t)\mid 0<s<t\}$. Another application
of the Dominated Convergence Theorem concludes the proof.\end{proof}

For the next result we need the following subspace of $\D([0,\mathscr{T}];\R)$:
\begin{align*}
\Fe([0,\mathscr{T}];(0,\infty)=\D([0,\mathscr{T}];\R)\cap\{f:[0,\mathscr{T}]\to(0,\infty)\mid f\text{ is non-decreasing}\}.
\end{align*}

\begin{lemma} \label{lemmacont}Let $g\in C((0,\infty)^{n};\R^{m}))$,
$m,n\in\N$ and $\mathscr{T}>0$. Then the mapping 
\begin{align*}
	\Psi & :(\Fe([0,\mathscr{T}];\R)^{n},d_{\infty}^{n})\to(\D([0,\mathscr{T}]^{n};\R^{m}),d_{\infty})\\
	x & \mapsto g(x_{1}(t_{1}),x_{2}(t_{2}),\dots,x_{n}(t_{n})),\quad\forall t=(t_{1},t_{2},\dots,t_{n})\in[0,\mathscr{T}]^{n},
\end{align*}
is continuous, where $d_{\infty}$ denotes the supremum metric.\end{lemma}

\begin{proof} We must show, that for every $x=(x_{1},\ldots,x_{n})\in\Fe([0,\mathscr{T}];\R)^{n}$
and $\varepsilon>0$ there exists a $\delta>0$ such that 
\begin{equation}
	d_{\infty}^{n}(x,y)<\delta\implies d_{\infty}(\Psi(x),\Psi(y))<\varepsilon.\label{contdef}
\end{equation}
To this end, note first that for all $i=1,\ldots,n$, $\inf_{t\in[0,\mathscr{T}]}x_{i}(t)=x_{i}(0)>0.$
Now, let $i=1,\ldots,n$ and consider the set $B_{x_{i}}=\{z\in\Fe([0,\mathscr{T}];\R)\mid d_{\infty}(x_{i},y)\leq x_{i}(0)/2\}$.
For $z\in B_{x_{i}}$ we have 
\[ |\sup_{t\in[0,\mathscr{T}]}\left|x_{i}(t)\right|-\sup_{t\in[0,\mathscr{T}]}\left|z(t)\right|| =\sup_{t\in[0,\mathscr{T}]}\left|x_{i}(t)-z(t)\right|<\max_{i}x_{i}(0)/2=:c_{x_{i}} \]
This implies $\sup_{t\in[0,\mathscr{T}]}\left|z(t)\right|\leq\sup_{t\in[0,\mathscr{T}]}\left|x_{i}(t)\right|+c_{x_{i}}<\infty$. Now set $C_{x_{i}}=c_{x_{i}}+\sup_{t\in[0,\mathscr{T}]}x_{i}(t)$.
Then, $z(t),x_{i}(t)\in[c_{x_{i}},C_{x_{i}}]$ for every $t\in [0,\mathscr{T}]$. Let
$K=\prod_{i=1}^{n}[c_{x_{i}},C_{x_{i}}]\subset(0,\infty)^{n}$. By
the Heine-Cantor Theorem, it follows that the restriction of $g$
to $K$ is uniformly continuous. This means that we can find $\delta_{0}>0$
such that for all $|t-s|<\delta_{0}\implies|g(t)-g(s)|<\varepsilon$.
Now, put $\delta=\min(\delta_{0},c_{x_{1}},\ldots,c_{x_{n}})$. We
conclude from above that if $y=(y_{1},\ldots,y_{n})\in\Fe([0,\mathscr{T}];\R)^{n}$
and $d_{\infty}^{n}(x,y)<\delta$ then $x(t),y(t)\in K$ for all $0\leq t\leq\mathscr{T}$,
from which (\ref{contdef}) follows. \end{proof}

\begin{proof}[Proof of Theorem \ref{consistencythm}] Let $\mathscr{T}>t_{0}$.
From Lemma \ref{lemma:limits}, the mapping $(s,t)\mapsto\psi(s,t;u,v)$
is continuous in $(0,\infty)^{2}$. Thus, we deduce that $\psi$ extends
(as in Lemma \ref{lemmacont}) to a continuous function $\Psi:\Fe([0,\mathscr{T}-t_{0}];(0,\infty)^{2})\to\D([0,\mathscr{T}]^{2};\R)$.
Moreover
\[
C^{n}(\cdot,\cdot;u,v)=\Psi((S_{t_{0}}[X]^n,S_{t_{0}}[X]^n)));\,\,\,C(\cdot,\cdot;u,v)=\Psi(S_{t_{0}}T,S_{t_{0}}T),
\]
where $S_{t_{0}}:\D([0,\mathscr{T}];\R)\to\D([0,\mathscr{T}-t_{0}];\R)$
denotes the shift operator, i.e.
\begin{align*}
	x & \mapsto(S_{t_{0}}x)(t)=x(t+t_{0}),\quad\forall t\in[0,\mathscr{T}-t_{0}].
\end{align*}
In view that $S_{t_{0}}$ is a continuous operator from $\D([0,\mathscr{T}];\R)$
to $\D([0,\mathscr{T}-t_{0}];\R)$ and $[X]^n\overset{u.c.p}{\Rightarrow}T$,
as $n\rightarrow\infty$, we can now apply The Continuous Mapping
Theorem (see for instance \cite{Cherubinietal12}) to conclude that
\[
d_{\infty}(\Psi(S_{t_{0}}[X]^n,S_{t_{0}}[X]^n),\Psi(S_{t_{0}}T,S_{t_{0}}T))\overset{\mathbb{P}}{\rightarrow}0,
\]
which is exactly the conclusion of the theorem.\end{proof}

\begin{proof}[Proof of Theorem \ref{clthm}]First note that thanks
to Lemma 5.3.12 in \cite{jacod2011discretization} we may and do
assume that $\left|\sigma_{t}\right|\leq C$ for some deterministic
constant $C>0$. Now, let 
\[
Z_{n}:=\left([X]^n_{t},[X]^n_{s}\right);\,\,\,Z:=\left(T_{t},T_{s}\right).
\]
Since $\mathbb{P}(\sigma_{t}^{2}>0)=1$ for all $t\geq0,$ we can
find $n$ large enough such that $Z_{n}\in \{(s,t):0<s,t,t\neq s\}$. Moreover, by Taylor's Theorem 
\[
\sqrt{n}\left[C^{n}(s,t;u,v)-C(s,t;u,v),\right]=\int_{0}^{1}\nabla\psi(Z+y(Z_{n}-Z),u,v)dy\cdot\sqrt{n}(Z_{n}-Z).
\]
From Theorem 5.4.2 in \cite{jacod2011discretization}, it follows
that for any $s\neq t$ 
\[
\sqrt{n}(Z_{n}-Z)\xrightarrow{s.d}\sqrt{2}\left(\int_{0}^{t}\sigma_{r}^{2}dW^\prime_{r},\int_{0}^{s}\sigma_{r}^{2}dW^\prime_{r}\right),
\]
where $W^\prime$ is a Brownian motion independent of $\mathcal{F}$. Therefore,
it is enough to show that 
\begin{equation}
	\int_{0}^{1}\nabla\psi(Z+y(Z_{n}-Z),u,v)dy\xrightarrow{\mathbb{P}}\nabla\psi\left(T_{t},T_{s},u,v\right).\label{eq:convergenderi}
\end{equation}
In view that $Z_{n}\overset{\mathbb{P}}{\rightarrow}Z$, every subsequence
$Z_{n_{k}}$ contains a further subsequence $Z_{n_{k(i)}}$ such that
$Z_{n_{k(i)}}\xrightarrow{a.s.}Z$. Fix $\omega\in\Omega_{t,s}:=\{\omega\in\Omega:Z_{n_{k(i)}}(\omega)\rightarrow Z(\omega),Z(\omega)\in \{(s,t)\mid 0<s<t\}\}.$
By using that $Z(\omega)\in V$, we can find an open ball with center
$Z(\omega)$ and radius $\rho(\omega)>0$ which is totally contained
in $\{(s,t)\mid 0<s<t\}$. Moreover, for every $\rho>\varepsilon>0$ there is $n_{0}\equiv n_{0}(\omega)\in\mathbb{N}$
such that 
\[
\left\Vert Z_{n_{k(i)}}(\omega)-Z(\omega)\right\Vert <\varepsilon,\,\,\,\forall\,n_{k(i)}\geq n_{0}.
\]
This in particular implies that for all $0\leq y\leq1$ and $n_{k(i)}\geq n_{0}$,
$Z(\omega)+y(Z_{n_{k(i)}}(\omega)-Z(\omega))$ is contained in an
open ball with center $Z(\omega)$ and radius $\varepsilon>0$. Therefore,
we can find a compact set $K_{\omega}\subseteq \{(s,t)\mid 0<s<t\}$ such that $Z(\omega)+y(Z_{n_{k(i)}}(\omega)-Z(\omega))\in K_{\omega}$
for all $0\leq y\leq1$ and $n_{k(i)}\geq n_{0}$. Hence,
by the continuity of $\nabla\psi$ on $\{(s,t)\mid 0<s<t\}$ (see Lemma \ref{lemma:limits})
and the Dominated Convergence Theorem, we deduce that as $n_{k(i)}\rightarrow\infty$
\[
\int_{0}^{1}\nabla\psi(Z(\omega)+y(Z_{n_{k(i)}}(\omega)-Z(\omega));u,v)dy\rightarrow\nabla\psi(Z(\omega);u,v),\,\,\forall\omega\in\Omega_{t,s}.
\]
(\ref{eq:convergenderi}) follows now by Theorem 6.3.1 in \cite{resnick2019probability}.\end{proof}

\bibliographystyle{plain}
\bibliography{ref}

\end{document}